\documentclass[11pt, a4paper, reqno]{amsart}

\usepackage{amssymb}
\usepackage{mathrsfs}

\usepackage[breaklinks]{hyperref}

\hypersetup{
unicode=true,
colorlinks=true,
linkcolor=blue,
citecolor=blue,
urlcolor=blue,
filecolor=blue,
bookmarksnumbered=true,
pdfstartview=FitH,
pdfhighlight=/N
}

\theoremstyle{plain}

\newtheorem{thm}{Theorem}

\newtheorem{conj}{Conjecture}

\theoremstyle{definition}

\newcommand{\ol}{\overline}

\def\({\left(}
\def\){\right)}
\def\Re{\operatorname{Re}}
\def\mcap{\operatorname{cap}}
\def\cop{\operatorname{cap}}
\def\const{\operatorname{const}}
\def\supp{\operatorname{supp}}

\def\CC{\mathbb C}
\def\PP{\mathbb P}
\def\MM{\mathrm M}

\let\leq\leqslant
\let\geq\geqslant
\let\le\leqslant
\let\ge\geqslant

\let\myo\overline
\let\mc\mathcal
\let\mbb\mathbb

\begin{document}

\title{Gonchar--Stahl's $\rho^2$-theorem and associated directions in the theory of
rational approximation of analytic functions}

\author{E.\,A.~Rakhmanov}
\address{University of South Florida, USA;
Steklov Mathematical Institute of Russian Academy of Sciences, Russia}
\email{rakhmano@mail.usf.edu}
\thanks{The work is supported by the Russian Science Foundation (RSF)
under a grant 14-50-00005.}

\date{23.03.2015}

\begin{abstract}
Gonchar--Stahl's $\rho^2$-theorem characterizes the rate of convergence of best
uniform (Chebyshev) rational approximations (with free poles) for one basic
class of analytic functions. The theorem itself, its modifications and generalizations,
methods involved in the proof and other related details constitute an important subfield
in the theory of rational approximations of analytic functions and complex analysis.

The paper briefly outlines essentials of the subfield. Fundamental contributions by
A.~A.~Gonchar and H.~Stahl are in the center of the exposition.

Bibliography:~\cite{Wals61Rus} items.
\end{abstract}

\maketitle


\markright{GONCHAR--STAHL $\rho^2$-THEOREM}

\setcounter{tocdepth}{2}
\tableofcontents


\section{Introduction. Statement of the theorem}\label{s1}

Gonchar--Stahl theorem characterizes the rate of convergence of best
uniform rational approximations (with free poles) for one basic class
of analytic functions. Its proof combines constructions and methods
from different branches of classical analysis and approximation theory.
This variety of significant connections explains the fundamental role
of the theorem in approximation theory and complex analysis. Some
of those facts and connections are briefly discussed below.

\subsection{Walsh theorem}\label{s1s1}
One of the main predecessors of Gonchar--Stahl theorem
is the well known J.~L.~Walsh theorem~\cite{Wals34} of the 1930s on the approximation of an
arbitrary element $f$ of an analytic function on a continuum $E$ of the (extended)
complex plane (see also the book~\cite{Wals60} and~\cite{Gon68}).

We consider the distance from $f$ to the class of rational functions of order $n$ in the
uniform metric on $E$,
\begin{equation}
\label{Rho}
\rho_n(f)= \rho_n(f, E)=\min_{r\in\mbb R_n}\max_{z\in E}|f(z)-r(z)|,
\end{equation}
where $\mbb{R}_n$ is the set of all rational functions $r_n=P_n/Q_n$ of order $\le n$ ($P_n,Q_n\in\mbb{P}_n $ -- polynomials of degree $\leq n$). It was a well known fact that $\rho_n(f) \to 0$ as $n\to\infty$. The problem was to determine the rate of convergence. Further, it was known that the rate is at least geometric and so, the problem was to estimate (asymptotically) the order of the associated progression.

The assumption that $f$ is analytic on $E$ means that there exists a domain
$\Omega$ containing $E$ such that $f\in H(\Omega)$; that is, $f$ is
holomorphic (analytic and single-valued) in $\Omega$.

Walsh's theorem asserts that $f\in H(\Omega)$ implies that
\begin{equation}
\label{LimR}
\varlimsup_{n\to\infty}
\rho_n(f)^{1/n} \le\rho= e^{-1/C(E,F)},
\end{equation}
where $F = \partial \Omega$ is the boundary of $\Omega$ and $C(E,F)$ is
the capacity of the condenser $(E,F)$. We will assume further that the
continuum $E$ has connected complement. If in addition $\Omega$ is
simply connected, then the number $1/\rho = 1/\rho(E,F)$ is also known as the
modulus of the ring domain $\Omega \setminus  E$ (there is a conformal
mapping of $\Omega$ on $\{\rho<|z|<1\}$).

In 1959 V.~D.~Erohin~\cite{Ero59} presented some examples proving that
this estimate is sharp; that is, it cannot be improved without further restrictions on $f$.
In particular, he constructed a function $f\in H(\Omega)$,
where $\Omega= \{z:|z|<R\}$ is a disk, such that for approximations of
$f$ on a smaller disk $E= \overline D_r = \{z:|z| \leq r\}$ equality in
\eqref{Rho} holds (it follows that for this function $f$ rational
approximations are not essentially better than polynomial
approximations of the same degree). This construction may be modified to
prove that for any domain $\Omega$ and any continuum $E\subset \Omega$
there exist $ f\in H(\Omega)$, such that equality in \eqref{Rho} holds (for
further details see the book~\cite{Wals60} by Walsh and a review by
S.~N.~Mergelyan included in the Russian translation of this
book~\cite{Wals61Rus} concerning the progress in this direction made
approximately between 1930 and 1960).

We note that in all those example there is only a small subsequence of
natural numbers along which the equality in \eqref{LimR} is reached. It was
determined later that ``in average'' the rate of convergence is essentially
better than in \eqref{LimR} for any function $f\in H(\Omega)$; see
Section~\ref{s1s2} below.

This means that the whole class $H(\Omega)$ always contains functions
with significantly irregular behavior of the sequence of rational
approximations. In this context Gonchar--Stahl's theorem essentially
asserts that the behavior of rational approximations is necessarily
regular for functions which have unlimited analytic continuation
outside a set of singularities of zero capacity. To formally state the
theorem we need the following definitions.

\subsection{Gonchar--Stahl \texorpdfstring{$\rho^2$}{rho2}-theorem}\label{s1s2}

Let $f|_E$ be an element of analytic function, which we want to
approximate. If analytic continuation of this element outside of $E$
has branch points, then there exist different (maximal) domains $\Omega$,
where $f|_E$ has holomorphic extension. Any of those domains may be
used in combination with the Walsh inequality \eqref{LimR} and,
therefore, $\rho$ in this inequality may be replaced by the following
constant $\rho(f)$ called {\it extremal modulus} associated with the
holomorphic function (analytic element) $f \in H(E)$
\begin{equation}
\label{ExtRho}
\rho(f)= \rho(f, E)=\inf \{e^{-1/C(E,F)}: F = \partial \Omega,\, f\in H(\Omega),\, E \subset \Omega \subset \mathcal D \}.
\end{equation}
There also exists a unique {\it extremal domain} $\Omega$, which
satisfies the condition $\rho(f, E)=\inf \{\exp\{-1/C(E,\partial\Omega)\}\}$
and is maximal among all such domains with this condition.

Let $\mathcal A (\mathcal D)$ be the class of all analytic elements in a domain $\mathcal D \subset \overline \CC$,
which admit analytic continuation along any path in this domain. We note that
this is one of the basic classes of analytic functions. For instance, solutions of
differential equations with coefficients in $H(\mathcal D)$ belong to
$\mathcal A (\mathcal D)$ (things are similar for all natural classes of equations).

The following is the Gonchar--Stahl $\rho^2$-theorem

\begin{thm} \label{GS}
Let $\mathcal D $ be a domain in the extended plane such that $\mcap(\overline\CC
\setminus \mathcal D) = 0$. Let  $E\subset \mathcal D$ be a continuum with
connected complement and $f \in H(E) \cap \mathcal A (\mathcal D)$. Then
\begin{equation}
\label{RhoSqu}
\lim_{n\to\infty} \rho_n(f, E)^{1/n}= \rho(f, E)^2
\end{equation}
In particular, the limit in the left hand side exists.
\end{thm}

This is not the most general form of the theorem; the conditions on $E$ may be
essentially relaxed, but for the purposes of our discussion this version is
representative enough. In the sequel we consider mostly cases when $E$ is an
interval or a disk (including the degenerated or local case when $E$ is a point).

There is also Gonchar's earlier version~\cite{Gon78} of the theorem
related to Markov-type functions. The Markov case is simpler, but it is
not a particular case of the Theorem~\ref{GS} above (see
also~\cite{GoLo78}). This case is discussed in some details separately
in Section~\ref{Sec 2}.

The proof of the general version of the theorem was presented in papers
by H.~Stahl~\cite{Sta86b} and A.~A.~Gonchar and the author~\cite{GoRa87}.
In~\cite{Sta86b} the upper estimate in~\eqref{RhoSqu} was proven and
in~\cite{GoRa87}  the corresponding lower estimate was obtained.
Actually, a general method of rational approximation was developed
in those two papers and solutions of several longstanding problems in
approximation theory were obtained there as an immediate
application of the method. The $\rho^2$-theorem was one of the classical
applications of the method. Stahl's paper~\cite{Sta86b}
contains another classical result; namely,  Stahl's celebrated theorem on the
convergence of diagonal Pad\'e approximants for functions with branch
points (corresponding with the degenerate case $E = \{\infty\}$ in the
$\rho^2$-theorem). The so-called ``$1/9$ problem'' on best
rational approximation of the exponential function on a semi-axis was
solved in~\cite{GoRa87} (it was obtained as a corollary of a general
theorem on Chebyshev approximation of a sequence of analytic
functions). Many other problems were investigated later by the same
method or its modifications (see, for example,~\cite{DeHuKu10},~\cite{Gon03}).
Further applications and generalizations
are coming (some of them are discussed below in this paper).                                        
For details and other references see the recent
reviews~\cite{Gon03},~\cite{ApBuMaSu11},~\cite{Rak12}, and also
JAT~\cite{Gon13}, \cite{Sta14}.

Our purpose is to present the essentials of the method mentioned
above. It is not possible to cover all the significant details. So,
we have selected some of them for the discussion. To this end we use the
$\rho^2$-theorem and its version for Markov-type functions as a starting
point of this discussion. It is also used to create a general context in the
first part of the paper.

\subsection{Contents of the paper}
In the next Section~\ref{Sec 2} we outline the proof of Gonchar's
$\rho^2$-theorem for Markov-type functions (Theorem~\ref{MarG} below).
The proof of this theorem \cite{Gon78} is essentially simpler than the proof of
Theorem~\ref{GS}. The simplicity of this situation makes it possible to
briefly discuss the main components of the proof in some details. In this
connection we also mention Stahl's theorem on the rate of rational
approximation of $|x|$ on $[-1,1]$ and results and conjectures by
Gonchar on the problem of characterizing classes of analytic
functions by the rate of their rational approximation.

In Section~\ref{s3} we discuss some details related to Stahl's theorem on the convergence of Pad\'e approximants for functions with branch points and, in particular, the characterization of minimal compact sets for functions with a finite number of branch points. This is a good ground for illustrating the geometric component of the method.
Some new results are then presented on stability of convergence in Stahl's theorem.
Finally some new conjectures are presented on asymptotics of complex orthogonal
polynomials (related to problems of convergence of Hermite--Pad\'e
approximations).

In the rest of this introduction we make some comments on the main
components of the method of proof of the $\rho^2$-theorem. Hopefully this
brief review may  represent to some extent the essence of the method. In the
conclusion we also make a few more comments on its connections and
applications.

\subsection{Brief description of the method}

Now, we describe briefly the main components of the method including, in
particular, a construction of ``near best'' approximations.

\subsubsection{Interpolation by rational function with free poles}
We begin by selecting a triangular table of interpolation nodes
on $E$, whose $n$-th row contains $2n+1$ points, and find the corresponding
rational function of order $n$ interpolating $f$ at the selected nodes. Now
we need to select nodes to obtain ``near best'' approximation.

We note that the proof of Walsh's theorem was based on interpolation with
fixed poles. In our ``free pole approximations'' method, poles are not known
in advance and have to be determined from the interpolating table. To decide
which nodes are ``near optimal'' (see~\cite{Gon03}) we need precise information on the location of
the poles of the interpolating functions.

The fundamental fact is that the denominators of the approximations
(polynomials whose roots are poles of the approximations) satisfy certain
complex (nonhermitian) orthogonality relations. Information on the poles
has to be derived from those orthogonality relations and we come to a
problem of asymptotics (in a typical case is the weak-$*$ zero
distribution) for a sequence of complex orthogonal polynomials.

\subsubsection{Asymptotics of complex orthogonal polynomials}
\ The method for studying the zero distribution of complex orthogonal polynomials
based on ingenious potential theoretic arguments was created by H.~Stahl
in~\cite{Sta86a},~\cite{Sta86b} and then substantially generalized in~\cite{GoRa87}.
The generalized Stahl method (the so-called GRS-method)
reduces the problem of asymptotics of orthogonal polynomials to an
equilibrium problem for the logarithmic potential.

This equilibrium problem is essentially different from ``standard''
equilibrium problems of complex analysis related to minimization of an energy
functional in the class of measures on a compact ``conductor''. Robin's measure of a
compact set in the plane and the equilibrium (signed) measure (distribution) of a
conductor (a pair of disjoint compact sets) are classical examples. The proof of Walsh's
theorem was based on condenser equilibrium distribution.

In the case of complex orthogonal polynomials we encounter a different kind of equilibrium which may be defined as equilibrium in a conducting domain.

\subsubsection{Equilibrium in a conducting domain}
This is a class of problems which may be described as problems of
critical points of an energy functional. Any such critical point is
naturally interpreted as an equilibrium position of an ``electric
charge'' on an open ``conductor''; the associated equilibrium is always
unstable. A classical example of such problems is Chebotarev's
problem, which asks for a continuum of minimal logarithmic capacity in the
class of all continua containing a given set of points. The Robin measure
of Chebotarev's continuum is an equilibrium distribution in the conducting
plane with a finite number of insulating points.

In particular, the proof of the $\rho^2$-theorem is related to some
generalization of the Chebotarev's problem (for the Green potential). A local
version of the problem (Pad\'e approximants) is related to a
generalization of the standard Chebotarev problem. The solution of
Chebotarev-type problems may be normally described in terms of
trajectories of a quadratic differential associated with the problem.
In many cases there are equivalent reformulations in terms of moduli of
families of curves so that this component of the method belongs also to
geometric function theory or even to differential geometry.

In~\cite{Sta85a}--\cite{Sta86b} H.~Stahl observed that the potentials
of Chebotarev-type systems of curves satisfy certain symmetry property
(now called $S$-property, i.e., the equality of the normal derivatives of the
equilibrium potential in two opposite directions) and this property may
be directly used to study complex orthogonal polynomials. He
generalized the problem to include curves with $S$-properties
($S$-curves) related to extremal cuts for Pad\'e approximants of
functions with branch points.

In~\cite{GoRa87} $S$-curves with an external field were part of the method. In the
general situation, the existence of an $S$-curve in a given class of functions
may be the key part of the whole problem (see~\cite{Rak12}).

\subsubsection{Lower bounds for approximations}
The first three parts of the method produce an ``optimal'' sequence of
rational approximations to $f$. As an immediate corollary it gives the
upper bound in \eqref{RhoSqu}. The corresponding lower bound is obtained
in~\cite{GoRa87} using special properties of the constructed sequence of
optimal approximations. The argument used is rather general and may be
stated as a separate theorem. The idea of the method was contained in one
of  Gonchar's earlier papers~\cite{Gon67},~\cite{Gon72}, see also the review~\cite{Gon03}.

\subsubsection{Applications and connections}
The method outlined above in connection with the $\rho^2$-theorem has
potentially a larger circle of applications. The two components in the
proof are especially important for applications, they are the asymptotics of
complex orthogonal polynomials and the related $S$-equilibrium problems.

Orthogonal polynomials are the key to a great variety of applications.
Together with the traditional ones in approximation theory, numerical
analysis, and spectral theory, many new applications were found, in particular in
mathematical physics, in the last two or three decades.

New classes of problems enter the theory related with random matrices and
statistics among other fields. New methods have been created in the
theory of orthogonal polynomials, in particular, steepest descent for
matrix Riemann--Hilbert problems. New versions of old methods such as the
Liouville--Green--Steklov method (a.k.a. WKB) were developed. In all
these cases an $S$-equilibrium configuration presenting a geometrical
component of the problem plays a key role.

The geometrical component originated by the existence problem for
$S$-equilibrium configurations leads to an environment somewhat similar
to the one related to general moduli problems in geometric function
theory (moduli of families of curves, quadratic differentials, critical
trajectories). This part is often present in difficult problems. For
instance, many important questions on matrix $S$-problems related to the
study of Hermite--Pad\'e polynomials are open;
see~\cite{Apt08},~\cite{Rak12}.

Finally, equilibrium problems create certain connection to integrable
systems (solutions theory). In many cases such problems related to
approximations are similar to those that come from mathematical
physics. For instance, some equilibrium problems associated with the
$\rho^2$-theorem and its generalizations are surprisingly close to
problems originated in the study of KdV or NLS equations by means of the
inverse scattering transform method. For some further comments
see~\cite{Rak12}.

\section{Approximation of Markov-type functions}\label{Sec 2}

Markov-type functions $f(z)$, which we call $\MM$-functions in the sequel, are
Cauchy transforms of positive measures with compact support on the real line
\begin{equation}
\label{Mar}
f(z):=\int_F\frac{d\sigma(t)}{t-z},\qquad z\in \Omega = \overline\CC \setminus F ,
\end{equation}
where $F$ is a (finite) interval (we may consider that $F$ is the minimal
interval containing the support of $\sigma$).

In particular, an $\MM$-function $f$ belongs to $\mathcal A(\overline\CC \setminus  e), e\subset
\mbb R$, if $e$ is
finite and the jumps of $f$ across the branch-cuts on the real line have constant
argument and are integrable. Note that the branch-cuts along $\mbb R$
constitute the boundary of the associated extremal domain. Thus, the classes of
$\MM$-functions and $\mathcal A(\overline\CC \setminus  e)$ with $\cop e =0$ are
overlapping but none of them contains the other.

Now, we fix an interval $E$ of the real line not intersecting $F$ and consider
best rational approximations to $f$ on $E$. Let $\rho_n(f, E)$ be the
distance from $f$ to the class $\mbb R_n$ in the uniform norm on $E$ (see
\eqref{Rho}).

\subsection{Gonchar's \texorpdfstring{$\rho^2$}{rho2}-theorem for Markov functions}\label{subsec 2.1}
The following is Gonchar's version of the $\rho^2$-theorem for $\MM$-functions.

\begin{thm}\label{MarG}
Let $\sigma'(x) =d\sigma/dx>0$ almost everywhere on $F$, then
\begin{equation}
\label{RhoSqu_2}
\lim_{n\to\infty}\rho_n(f, E)^{1/n}= \rho(f, E)^2
\end{equation}
\end{thm}

Main components of the proof.

\subsubsection{Interpolation}
We begin with an arbitrary triangular table of points $\{\zeta_{k,n}\}
\subset E$, where $n$ is a natural number and $k = 1,2,\dots 2n$ for a fixed
$n$. We define $W_n(z) = \prod_{k=1}^{2n} (z-\zeta_{k,n})$.

Next we define the $n$-th order multipoint Pad\'e approximant $r_n$ to
$f$ associated with the interpolation table $\{\zeta_{k,n}\}$
(see, for example,~\cite{Gon75},~\cite{Gon75b},~\cite{Gon81}). For
technical reasons it is convenient to use exactly $2n$ interpolation points
and then use the interpolating function with the condition $r_n(\infty) =
0$.

For each $n$ there exist a pair of polynomials $P_n \in\PP_{n-1}$ and
$Q_n \in\PP_n$ such that $Q_n \not\equiv 0$ and the condition
$$
F_n(z) = (Q_n f - P_n)(z) / W_n (z) \  \in \ H(E)
$$
($F_n$ is analytic on $E$) is satisfied. Indeed, the last condition is equivalent to a system
of $2n$ linear homogeneous equations for $2n+1$ coefficients of the polynomials
$P_n \in\PP_{n-1}$ and $ Q_n \in\PP_n$ (in case of distinct nodes those equations are
 $(Q_n f - P_n)(\zeta_{k,n})= 0$, \ $k =1,2,\dots 2n$).
Such a system always has a nontrivial solution. This proves the existence of required
polynomials $P_n, \  Q_n \not\equiv 0$.

We set $r_n = P_n/Q_n$. This function does not necessarily interpolate $f$ at all
nodes (common zeros between $P_n$ and $Q_n$ are possible where the interpolation
may be lost). For Markov-type functions this cannot happen as it follows
from subsequent considerations.

\subsubsection{Orthogonality. Hermite interpolation formula}
The denominator $Q_n$ satisfies the following orthogonality conditions
\begin{equation}
\label{Ort}
\int_F Q_n(x)\, x^j \,\frac {d\sigma(t)}{W_n(t)} = 0 , \qquad j = 0,1, \dots, n-1.
\end{equation}
The following identity, called Hermite interpolation formula, is also important
\begin{equation}
\label{Her}
 f(z) - r_n(z) = \frac {W_n(z)}{Q^2_n(z)}\int_F \, \frac {Q^2_n(t)}{W_n(t)}\,\frac {d\sigma(t)}{z-t}= 0, \qquad z\in\Omega.
\end{equation}

The proof of \eqref{Ort} is obtained integrating $z^j\,F_n(z)$ over a
contour $C$ separating $F$ from $E$ and $\infty$. In particular,
\eqref{Ort} implies that the zeros of $Q_n$ are simple and belong to $F$. After that
\eqref{Her} follows from   Cauchy's integral representation for the
function $Q_nF\in H(\text{Ext}\,C)$, where $\text{Ext}\,C$ denotes
the unbounded connected component of the complement of the contour $C$.

\subsubsection{Zero distribution. Balayage}
Conditions \eqref{Ort} present a model
situation of orthogonality with varying weights on the real line. The assertion
on the asymptotics of the associated orthogonal polynomials $Q_n$ is formulated in
terms of weak-$*$ convergence of the normalized counting measure. The counting measure of a
polynomial $P$ is defined as the sum of unit masses at the zeros of $P$ (counting multiplicities) and is denoted
by $\mc{X}\(P\)=\sum_{P(\zeta)=0}\delta(\zeta)$ \ ($\delta(\zeta)$ is
the unit measure supported at the point $\zeta$).

Now we assume that the interpolation table has a limit distribution (limit
density) represented by a unit positive measure $\mu$ on $E$. More precisely,
this assumption means that the sequence $\mc{X}(W_n)/2n$ is weak-$*$ convergent
to the measure $\mu$ as $n \to \infty$. We denote this fact by $\frac1{2n}\
\mc{X}\(W_n\){\overset{*}{\to}}\mu$.

The basic fact regarding free poles real interpolation of $\MM$-functions is
that, if the interpolation table has limit density $\mu$, then the denominators
$Q_n$ have the limit distribution $\lambda$, which is the balayage of $\mu$
from $E$ onto $F$. Formally,
\begin{equation}
\label{Dist}
\frac1{2n}\ \mc{X}\(W_n\)\ \ {\overset{*}{\to}}\ \ \mu
\qquad \text{implies that}\qquad
\frac1{n}\ \mc{X}\(Q_n\)\ \ {\overset{*}{\to}}\ \ \lambda,
\end{equation}
where $\lambda$ is the unit measure on $F$ defined by the condition
\begin{equation}
\label{Bal}
U^\lambda(x) - U^\mu(x) = C_F = \const, \qquad x\in F
\end{equation}
(by $U^\nu(x) = - \int\log|x-t| d\nu(t)$ we denote the logarithmic potential of a measure $\nu$).

\subsubsection{Convergence. Upper bound for the rate on $E$}
Comparing the boundary values on $F$ and
singularities on $E$, it is directly verified that $U^{\mu - \lambda} (z) + C_F = G^\mu(z) =
\int g(z,t) d\mu(t)$ is Green's potential of $\mu$ with respect to the
domain $\Omega = \overline\CC \setminus  F$ \ ($g(x,t)$ is Green's function for
$\Omega$ with pole at $t$). It follows by \eqref{Her} and
\eqref{Dist} that the interpolating sequence $r_n(z)$ associated
with an interpolation table having limit density $\mu$ verifies
\begin{equation}
\label{Asy}
|f(z) - r_n(z)|^{1/n}\, \to \, \exp{\{- 2 G^\mu(z)}\}, \qquad
\text{uniformly for}\quad z\in \Omega \setminus  E
\end{equation}
and, moreover,
\begin{equation} \label{Asy-1}
\max_{x\in E}|f(x) - r_n(x)|^{1/n}\, \to \, \exp{\{- 2 \min_{x\in E} G^\mu(x)}\}.
\end{equation}
To obtain the best possible estimate from \eqref{Asy-1}, we need to find a measure $\mu$, which maximizes
$w(\mu) = \min_{x\in E} G^\mu (x)$
in the class of all unit measures $\mu$ on $E$. The problem is well known in classical complex analysis; its solution is $\mu = \lambda_E$, Green's equilibrium measure on $E$ relative to $\Omega = \overline\CC \setminus  F$. A characteristic property of this measure is
\begin{equation} \label{Bal-1}
U^\lambda(x) - U^\mu(x) = C_E = \const, \qquad x\in E .
\end{equation}
Relations \eqref{Bal} and \eqref{Bal-1} together mean that the pair of measures
$\lambda_E = \lambda$ and $\mu = \lambda_F$ form the equilibrium distribution for
the condenser $(E,F)$. The capacity $C(E,F)$ of the condenser is defined as
$C(E,F) = 1/w$. In terms of equilibrium constants in \eqref{Bal}, \eqref{Bal-1}
we have $w = C_F -C_E$.
From here $\varlimsup_{n \to \infty}\rho_n(f)^{1/n}\leq \rho(f)^2$.

We note that the equilibrium problem related to this situation is the
standard equilibrium for a plane condenser (pair of disjoint compact sets) and,
thus, we do not seem to have a nonstandard $S$-equilibrium problem here.
The reason is that on the real axis any equilibrium potential is
automatically symmetric with respect to the real axis (normal derivatives of a
potential in two opposite directions are equal). So, the associated $S$-property
is actually following from the symmetry of the situation.

\subsubsection{Lower bound and strong asymptotics}
It follows directly from \eqref{Her} that $f - r_n$ is real
and has exactly $2n$ zeros on $E$. Thus, it makes $2n+1$ oscillations on
$E$, whose amplitudes are asymptotically estimated by \eqref{Asy}. These
estimates make it possible to use the classical Ch.~De~la~Vallee-Poussin
inequality~\cite{Val11} (instead of Gonchar's general complex argument mentioned
above) to obtain the corresponding lower bound.

Moreover, under the Szeg\H{o} condition~\cite{Sze62} on the measure $\sigma$ it is
possible to slightly modify the interpolation nodes (the measure $\mu$) in
such a way, that the difference $f-r_n$ is asymptotically
equioscillating. Then, the application of Ch.~De~la~Vallee-Poussin's
estimates will infer the strong asymptotics for the error of the best
approximations
\begin{equation} \label{Asy_2}
\rho_n(f) \, = \, \gamma (\sigma, E)\, \rho(f)^{2n} (1+\varepsilon_n) ,
\end{equation}
where $\varepsilon_n \to 0$ as $n\to\infty$ and $ \gamma (\sigma, E)$
is an explicit constant (the same is true for the error of interpolation  $\max_{x\in E}|f(x) - r_n(x)|$).

Details related to this and other similar results may be found in the book by
H.~Stahl and V.~Totik~\cite{StTo92} on general orthogonal polynomials. One example of a more complex (but still Markov) situation presented by Stahl's theorem on the rational approximation of $\sqrt{x}$ on $[0,1]$ is in  Section \ref{subsec 2.2} below.

For a rather general class of analytic functions (including functions with complex branch points)
A. Aptekarev~\cite{Apt02} proved a theorem
on the exact constants of approximation by rational functions of order \ $\leq n$.
In particular, he obtained a relation of type~\eqref{Asy_2} for the functions from this class.
His theorem yields the following formula of strong asymptotics for the error
$\rho_n=\rho_n(e^{-x})$ of the best uniform approximation on the semi-axis $[0,+\infty)$
to the function~$e^{-x}$ with rational functions of order \ $\leq n$:
\begin{equation}
\rho_n=2 v^{n+1/2}(1+o(1)),\qquad
n\to\infty,
\label{Apt_Asym}
\end{equation}
where $v$ is the so-called Halphen constant (see~\cite{GoRa87},~\cite{Mag87}).
Formula~\eqref{Apt_Asym} proves a conjecture of Magnus~\cite{Mag87}
on the exact constant of rational approximation of function~$e^{-x}$ on the semi-axis $[0,+\infty)$.

The method in \cite{Apt02} was based on a study of strong asymptotics for
complex orthogonal polynomials using steepest descent for matrix Riemann--Hilbert representation of such polynomials.
In this connection see also \cite{ApYa11}.

See also \cite{MaRaSu12}, where both steepest descent for matrix Riemann--Hilbert and WKB are used.

\subsubsection{Generalization. Equilibrium measure}
Formulae~\eqref{Dist}--\eqref{Bal} were proven in the original papers~\cite{GoLo78}
and~\cite{Gon78} only for the case $\mu = \lambda_E$, which was enough to
conclude the proof of Theorem~\ref{MarG}.

More general orthogonal polynomials $Q_n$, defined by
\begin{equation}\label{Ort_1}
\int_F Q_n(x)\, e^{-2n \varphi_n(x)}x^j \,d\sigma(t) = 0 , \qquad j = 0,1, \dots, n-1 ,
\end{equation}
were studied in the paper~\cite{GoRa84} by A.~A.~Gonchar and the author with the following result (simplified version)

\begin{thm}\label{Ort-2}
Let $\sigma'(x) =d\sigma/dx > 0 $ almost everywhere on $F$ and suppose that the sequence
$\varphi_n(x)$ converges to $\varphi(x)$ uniformly on $F$. Then $\frac1{n}\
\mc{X}\(Q_n\)\ \ {\overset{*}{\to}}\ \ \lambda$,
where $\lambda=\lambda_\varphi$ is the equilibrium measure of $F$ in the external field $\varphi (z)$,
which is a unit measure on $F$ defined by
\begin{equation}
\label{EqM}
U^\lambda(x) + \varphi(x) = C = \const, \quad x\in \supp\lambda, \qquad U^\lambda(x) + \varphi(x) \geq C ,\quad x\in F.
\end{equation}
\end{thm}

This theorem was, probably, the first general result on zero distribution of
orthogonal polynomials with varying weights.

\subsection{Stahl's theorem on approximation of \texorpdfstring{$|x|$}{|x|} on \texorpdfstring{$[-1,1]$}{[-1,1]}}\label{subsec 2.2}
The problem of estimates for $\rho_n = \rho_n(|x|, [-1,1])$ was introduced
by D.~J.~Newman in 1964~\cite{New64}, who proved that $e^{- c_1\sqrt n} \leq
\rho_n \leq e^{ - c_2\sqrt n}$ with some $c_1 \geq c_2 > 0$.

It is easy to see that $\rho_{2n} = \rho_n(\sqrt x, [0,1])$, so the problem is reduced
to the approximation of $\sqrt x$ on $[0,1]$.

Representing the function $\sqrt z$ in of the domain $\{|z| <2\} \setminus  (-2, 0]$
by the Cauchy integral and defining
\begin{equation}\label{Root}
f(z):= \frac 1\pi \int_{[-1,0]}\frac{\sqrt{-t}}{z - t}\, dt,\qquad z\in \Omega = \overline\CC \setminus  [-1, 0]
\end{equation}
we find that $g(x) = \sqrt x - f(x)$ is analytic on $[0,1]$ and, therefore, its
rational approximations converge to $g$ geometrically. This, together with
Newman's estimate, implies that $ \rho_n(\sqrt x, [0,1])/ \rho_n(f(x),
[0,1]) \to 1 $ as $n \to \infty$.

So the problem is reduced to the study of the best rational approximation of the
$\MM$-function $f$ on $E = [0,1]$. Basically, we can use the method described
in Section~\ref{Sec 2} above, but this method has to be modified.

The problem now is more difficult than the problems
discussed in Section~\ref{Sec 2}. The condenser $(E, F)$ associated with
the current situation $E = [0,1]$, $F=[-1, 0]$ is degenerated, since the plates
$E$ and $F$ have a common point, the equilibrium $\lambda_E - \lambda_F $
distribution for such condenser does not exist (collapsing situation).

Stahl used the condenser with a logarithmic external field on the plate
$F$ (which comes from term $\sqrt{-t}$ in \eqref{Root}).
The external field prevents the equilibrium distribution from collapsing and the
weighted equilibrium distribution may be used to define an optimal interpolating table.
Stahl was able to obtain strong asymptotics for associated orthogonal polynomials. Then he obtained strong asymptotics for the error of approximation. As a result he proved in~\cite{Sta92} the following remarkable theorem

\begin{thm}
For $\rho_n = \rho_n(|x|, [-1,1])$ we have
\begin{equation}
\label{Root-1}
\lim_{n\to\infty} \rho_n\, e^{ \pi\sqrt n} = 8.
\end{equation}
\end{thm}

The result was conjectured by R.~Varga~\cite{VaRuCa91} on the basis of
numerical experiment. The correct constant $c_ 1 = c_2 = \pi$ was
earlier found by N.~S.~Vyacheslavov~\cite{Vya74}.

\subsection{Some problems and conjectures by Gonchar related to the \texorpdfstring{$\rho^2$}{rho2}-theorem}\label{2.3}
A broader context related to the $\rho^2$-theorem is the general problem of characterizing classes of analytic functions through the rate of convergence of their best rational approximations.

The corresponding problem for polynomial approximations essentially admits a
general solution and the associated theory is well known. For rational
approximations the situation is more complicated. It is usually difficult to
find a criterium in terms of best rational approximation since direct and
inverse theorems are mostly far away from each other
(cf.~\cite{Gon82},~\cite{Gon84-2},~\cite{Gon84-3}).

A typical example is related to characterizing the class of functions with supergeometric rate of convergence of best rational approximations
\begin{equation} \label{Rho-S}
\lim_{n\to\infty}\rho_n(f, E)^{1/n}= 0 .
\end{equation}
A direct theorem by Ch.~Pommerenke~\cite{Pom73} asserts that, if $f\in H(\overline\CC \setminus  e)$ and
$\cop(e) = 0$, then \eqref{Rho-S} is valid for any $E \subset \overline\CC \setminus e$.
The inverse is not true. Basically, knowing \eqref{Rho-S} we cannot
say anything about the set of singularities of $f$.

On the other hand, Gonchar proved in~\cite{Gon72} that \eqref{Rho-S}
implies that $f$ is quasianalytic (there is a uniqueness theorem for
such functions similar to the one for analytic functions). He also
proved~\cite{Gon72} that the function $f$ is single-valued in all of its
Weierstrass domain if \eqref{Rho-S} is satisfied. There are more
theorems by Gonchar of that kind; see~\cite{Gon74},~\cite{Gon74b} for
details.

Soon after the $\rho^2$-theorem for $\MM$-functions was proven, Gonchar raised
the following general question: for what kind of functions the limit of
$\rho_n(f, E)^{1/n}$ as $n\to\infty$ exists and is positive? In other words,
which  functions have regular behavior of the sequence of their best rational
approximations. His basic idea was that all the ``natural'' functions are
regular and for any such function we have
\begin{equation}
\label{Rho-3}
\lim_{n\to\infty}\rho_n(f, E)^{1/n}=\rho(f, E)^2
\end{equation}
for any continuum $E$ in the domain of the function (see
\cite{Gon78},~\cite{Gon84-1},~\cite{Gon84b},~\cite{Gon86}). In other words, if $\lim \rho_n(f, E)^{1/n}$ exists
then it is equal to $\rho(f)^2 $. All subsequent
results seem to confirm the conjecture but it is not clear how it may be
proved.

Anyway, Theorem~\ref{MarG} means that Markov-type functions are regular
(\eqref{Rho-3} is satisfied under mild restrictions on the measure), which
was also an important argument in favor of Gonchar's
$\rho^2$-conjecture that any function $f\in \mathcal A(\myo \CC \setminus  e)$
with a finite set $e$ is a regular function (a stronger version is contained in the
Gonchar--Stahl theorem). In particular, he also conjectured that all
algebraic functions are regular and for any $E$ free of singularities
we have the stronger estimates
\begin{equation}
\label{Rho-Gon}
0 < C_1(f,E) \leq \frac{\rho_n(f, E)}{\rho(f, E)^{2n}}\leq C_2(f,E)
\end{equation}
(this Gonchar conjecture seems to be still generally open, even though within
the range of existing methods).

Another general conjecture by Gonchar was that, if a function has worse than
$\rho^2$-rate of best rational approximation and
$\varlimsup\limits_{n\to\infty} \rho_n(f, E)^{1/n} > \rho(f, E)^2$, then it is ``not
regular'' and there is another subsequence where the rate is better than ``normal''; that is,
$\varliminf\limits_{n\to\infty}\rho_n(f, E)^{1/n}<\rho(f, E)^2$.
In particular, his conjecture was that for any function $f$ we have
\begin{equation}
\label{Rho-4}
\varliminf_{n\to\infty}\rho_n(f, E)^{1/n} \leq \rho(f, E)^2 \qquad \text{for any}\quad f\in H(E).
\end{equation}
Later this conjecture was proved by O.~G.~Parfenov \cite{Par86} and
V.~A.~Prokhorov \cite{Pro93}--\cite{Pro94} who also obtained a stronger
inequality
\begin{equation}
\label{Rho-5}
\varlimsup_{n\to\infty}\,\prod_{k=1}^n \rho_k(f, E)^{1/k} \leq \rho(f, E)^2.
\end{equation}
The proofs of the theorems by Parfenov~\cite{Par86} and
Prokhorov~\cite{Pro93}--\cite{Pro94} were based on a combination of fixed poles of
interpolation and theorems of singular numbers of Hankel operators. This is
essentially another important direction in approximation theory in many ways
different from the one under consideration and we do not go into further details.

\section{Stahl's theorem on Pad\'e approximants}\label{s3}

So far the rate of convergence of best rational approximations was
discussed. Now we pass to the convergence properties of the
approximating functions; see \eqref{Asy} and \eqref{Asy-1} as
examples. The construction of near-best rational approximation $r_n(z)$ to
$f$ in the context of Theorem \ref{GS} may also be arranged in such a way,
that these functions converge to $f$ uniformly on compact subsets of
the whole extremal domain of analyticity of $f$.

The convergence problem for rational approximations is more convenient to
discuss for the case of (diagonal) Pad\'e approximants, the best local
rational approximants to a power series. It is also convenient to select
an interpolation point at infinity, so that all branch points of the function
are finite.

\subsection{Pad\'e approximants for functions with branch points} \label{3.1}
 Let
\begin{equation}
\label{f}
f(z)=\sum^\infty_{n=0}\frac{f_k}{z^k}
\end{equation}
be a function analytic at infinity. Pad\'e approximants $\pi_n(z)=(P_n/Q_n)(z)$
to $f$ are defined by the condition
\begin{equation}
\label{eq2.2}
R_n(z): =\(Q_nf-P_n\)(z)=O\left(1/z^{n+1}\right),\quad z\to\infty,      
\end{equation}
where $P_n,Q_n\in\mbb{P}_n$ and $Q_n \not\equiv 0$ (see~\cite{BaGr81},
for details).
In connection with best rational approximations we note that the $\pi_n(z)$ are
limits of the best approximations on $D_R = \{z: |z| \geq R\}$ as $R\to
\infty$~\cite{Wals64}. In this sense Pad\'e approximants are the local version of
best rational ones. Function $R_n$ is called the remainder.                                          

\subsubsection{Markov theorem}
An old classical convergence theorem proved by
A.~A.~Markov~\cite{Mar95} in 1895 asserts that, if $f(z)$ is an $\MM$-function \eqref{Mar}, then
the associated sequence $\pi_n(z)$ converges to $f$ uniformly on compact subsets of the complement of
$F$ \ (the minimal interval containing the support $S_\sigma$ of the measure
$\sigma$). Note that $f$ is holomorphic in the domain $\Omega = \overline\CC \setminus
S_\sigma$, which may be larger than $\overline\CC \setminus  F$. The functions $\pi_n$ may
have poles in this larger domain, but they still converge there in capacity.

The fact that Pad\'e denominators $Q_n$ are orthogonal polynomials with
respect to $\sigma$ had been discovered earlier in 1855 by
P.~L.~Chebyshev~\cite{Che55}.

\subsubsection{Nuttall's minimal capacity conjecture}
One of the main problems
in the theory of Pad\'e approximants in the period 1960--1970 was the
convergence problem for functions with branch points. If element
\eqref{f} at infinity represents a function $f\in
\mc{A}(\mbb{C}\setminus {e})$, where $e$ is, say, a finite set of branch
points, then Pad\'e approximations to $f$ may converge to $f$ only in a
domain where $f$ is single-valued. What is actually the domain of
convergence?

The first results on the convergence of Pad\'e approximants for functions
with some special type of branch points were obtained by J.~Nuttall who
also made the following conjecture (see~\cite{NuSi77},~\cite{Nut84}).
Let $f\in \mc{A}(\mbb{C}\setminus {e})$ where $e$ is a finite set and
\begin{equation}
\label{F}
\mathcal F =\{F\subset\mbb{C}:f\in H(\mbb{C}\setminus  F)\}
\end{equation}
be the set of compact cuts $F$, which makes $f$ single-valued.
Let, further, $F_f\in \mathcal F $ be the cut of minimal capacity
\begin{equation}
\label{F-f}
\cop(F_f)=\min_{F \in \mathcal F }\cop(F) .
\end{equation}
Nuttall's main conjecture was that the sequence $\left\{\pi_n\right\}$ converges to $f$ in capacity in the complement to $F_f$
$$
\pi_n\overset{\cop}{\to}f,\quad z\in\mbb{C}\setminus  F_f.
$$
He also formulated a conjecture on strong asymptotics of Pad\'e
denominators, which he proved in some particular
situations~\cite{NuSi77},~\cite{Nut84}.

\subsection{Stahl's theorem} \label{3.2}
A general theorem on the convergence of Pad\'e approximants including, in particular,
Nuttall's conjecture, was proven by
H.~Stahl~\cite{Sta85a}--\cite{Sta86a}. Here is the original statement of the
theorem where the compact set $F_f$ of minimal capacity is characterized
equivalently in terms of the $S$-property.

\begin{thm}
Let $e$ be a compact set of zero (logarithmic) capacity $ \cop e = 0$ and
$f\in \mathcal A_e=\mathcal A(\overline\CC\setminus  e) $ is not single-valued in
$\overline\CC\setminus  e$. Then the following assertions $(A)$, $(B)$ and $(C)$ hold

$(A)$ There exists a unique compact set $F = F_f$ in the plane, which is the
union of analytic arcs (up to subsets of capacity zero), with the
following properties. The complement of $F$ is connected, $f$ is single
valued in $\overline\CC\setminus  F$ (so that $F\in\mathcal F (f)$), the jump of
$f$ across any arc in $F$ is not identical zero and, finally, the equality
\begin{equation}
\label{S}
 \frac{\partial g}{\partial n_1}(z) = \frac{\partial g}{\partial n_2}(z), \quad z\in F^0
\end{equation}
(called $S$-property) holds for the Green function $g=g(z,\infty)$ of $\CC\setminus  F$ with pole at infinity,
where $F^0$ is the union of the open parts of the arcs
constituting $F$ ($n_1, n_2$ are two oppositely directed normals to
$F^0$ at the point~$z$).

$(B)$ For the Pad\'e denominators $Q_n$ associated with $f$ we have $\mc X(Q_n)/n
\xrightarrow{*} \lambda$, where $\lambda=\lambda_F$ is the Robin measure of
the compact set $F$.

$(C)$ The sequence of Pad\'e approximants $\pi_n = P_n/Q_n$ associated with
$f$ converges in capacity to the function $f$ inside (i.e., on compact
subsets) of the domain $D:=\overline\CC\setminus F$.
\end{thm}

The exact rate of convergence in capacity was also included in the theorem,
but our further discussion is related to assertions $(A)$ and $(B)$.
Part $(C)$ is essentially a corollary of $(B)$.

The most important part of the theorem is part $(B)$. Rather sophisticated and
entirely original potential theoretic methods were used in this part of the
proof. The starting point was the following orthogonality condition for Pad\'e
denominators $Q_n$
\begin{equation} \label{Ort_3}
\oint_F Q_n(z) z^k f(z) dz = 0, \quad k =0,1, \dots, n-1,
\end{equation}
where integration is taken over any system of contours separating $F$
from infinity. Very important is that Stahl's proof was the first instance
of the effective use of complex orthogonality and this was a
significant breakthrough in the theory.

Another interesting fact about Stahl's proof is that additional assumptions
do not lead to any simplifications. The proof of the theorem for a set
$e$ with three branch points (not on the line) is identical to the original
proof for sets $e$ of capacity zero. Additional assumptions on the
character of the branch points do not bring any simplifications either (one
exception is a square root of a rational function). It seems that this part
of the theorem does not have any simple complex particular cases for branch
points not on the line. The case when the branch points are on a line may be
essentially viewed as part of the Markov theorem (some additional
assumptions are formally needed).

It is not possible here to go into any further details related to this part of
the proof. We have in mind that the Robin measure $\lambda=\lambda_F$ of the
extremal compact set $F = F_f$ represents the limit zero distribution of the Pad\'e
denominators and concentrate on the further characterization of $F$ (geometric
component of the problem). Part $(A)$ of the theorem essentially defines $F$ by
the $S$-property.

Recall that the similar part in the proof of the $\rho^2$-theorem was
represented by a Green equilibrium problem. In more general situations, the
geometric component may be represented by a more general kind of $S$-equilibrium
problems. Then existence itself would be a problem.
Constructive solutions is another problem. To some extent,
further progress in the theory depends on the development of the geometric
component of the method.

\subsection{`Geometry' of Stahl's theorem} \label{3.3}
The geometric part in Stahl's theorem is a particular case of a general
$S$-equilibrium problem; it is the case of single logarithmic potential and the zero external field.
The extremal compact set $F_f$ in such settings always exist and it has
comparatively simple constructive characterization.
Similar characterization (in case of existence) may be obtained in
more general cases of pure logarithmic or Green potential and an external
field, which is harmonic outside of a set of capacity zero. Other cases are
essentially open; see~\cite{Rak12}.
Next we discuss a case when set $e$ of singularities of $f$ is finite (and there is no zero external field).
Formulae related to the case are explicit. At the same time the case is still representative and may shed light on the nature of the problem at large.

\subsubsection{Quadratic differential}\label{s3.3.1}
Let  $e = \{a_1, a_2, \dots, a_p\}$ be finite set of distinct points; we denote
$A(z) = (z-a_1)(z- a_2) \cdots(z - a_p)$. Let
$f\in \mathcal A_e =\mathcal A(\overline\CC \setminus  e)$ and $F_f$ be the associated extremal compact.  In the sequel we call it $S$-compact or Stahl compact. The following characterization
of $F_f$ is valid.

There exists a polynomial $V_f$
$$
V_f(z) = (z-v_1) (z-v_2) \dots (z-v_{p-2}) , \qquad \text {where}\quad v_j = v_j(f)
$$
depending on $f$ and $e$ of degree $p-2$ such that the $S$-compact set $F_f$ is the union of some critical trajectories of the quadratic differential $ - (V/A)\,(dz)^2$, where $V = V_f$.

The assertion follows from Stahl's results~\cite{Sta85a}--\cite{Sta85c}. However, the statement is close
to some traditional theorems in geometric function theory~\cite{Str84}.
An alternative proof based on a Max-Min energy problem was presented
in~\cite{PeRa94}; see the review~\cite{Rak12} for more details.

In addition, the following condition is satisfied: $ - (V/A)\,(dz)^2$ is
the quadratic differential with closed trajectories (in the terminology
of~\cite{Str84}). In our particular case it means that all its
trajectories, defined by the inequality $- ( V(z)/A(z)) (dz)^2 > 0 $ are either closed
contours or critical-analytic arcs, connecting some pair of zeros of
$AV$.

Moreover, function $\sqrt{ {V(z)}/{A(z)}}$ has a holomorphic branch in $\Omega = \overline\CC\setminus F_f$ and the Green function $g$ for $\Omega $ with pole at infinity can be written as
\begin{equation} \label{Gre}
g(z) = \Re G(z),\qquad G(z) = \int_a^z \, \sqrt{\frac {V(t)}{A(t)}}\, dt \quad (a\in e),
\end{equation}
(the branch of the root is such that $g(z) = \log|z| + o(1)$ at infinity).
The $S$-property \eqref{S} of the Green function follows directly from this representation.

Representation \eqref{Gre} establishes one-to-one correspondence between $S$-compacta $F_f$ and polynomials $V_f$. Zeros of  $V_f$ may serve as coordinates of  $F_f$.

\subsubsection{Family of polynomials \ $V_f$, \ $f \in \mathcal A_e$  }\label{3.3.2}
The problem of a constructive determination of the compact $F_f$ for a given
$f$ has two components. This compact  depends, first of all, on the branch set
$e$ of the function $f\in\mathcal A_e$. It depends also on the branch type of the function
(determined by indicating the loops, along which the analytic
continuation of the original element given at the infinity point leaves
it unchanged). It is convenient to separate the two dependencies by introducing the family of compacta $F_f$ associated with all functions $f \in \mathcal A_e$ having fixed set $e$ of branch points.

It is not difficult to prove that this family is finite. The number of its elements depends on the number of points in the set $e$ and their configuration (we do not discuss calculation of this number). Now we concentrate on characterization of this family.

Since each compact $F_f,$ is uniquely defined by associated polynomial $V_f$, the whole family
$F_f,$ \ $f \in \mathcal A_e$ may be described in terms of associated family of polynomials $V_f$, which we denote by
$$
\widetilde {V} (e) = \left\{ V_f : \,\, f \in \mathcal A_e \right\}
$$

Polynomial $V \in \widetilde {V} (e)$ is determined by its roots, that is, by $p-2$ complex numbers $v_j$
playing role of coordinates. We may ask, therefore, if some kind of equations may be written in terms of coordinates $v_j$.  Some equations may, indeed, be derived from the characterization of $S$-compacta as critical trajectories of quadratic differentials (see Section \ref{s3.3.1} above). Those equations (written in terms of periods of quadratic differentials) belong to a well known class of equations and such equations are usually not easy to deal with. In particular, there is a difficult combinatorial element in their structure and the detailed analysis of the situation may not be presented here. Below we outline briefly two possible two ways the problem may be approached without going into all the details.

In the next subsection, we introduce a family of hyperelliptic Riemann surfaces associated with the family $F_f,$ \ $f \in \mathcal A_e$ of Stahl compacta. In terms of this family of Riemann surfaces we define a mapping in the set of monic polynomials of degree $p-2$. Then polynomials $V\in\widetilde {V} (e)$ are defined as fixed points of this mapping. It seems that nothing related to this approach has been published so far.

It is possible that a natural way to generalize Stahl theorem for Hermite--Pad\'e approximation goes through a proper generalization of this approach. Anyway, the associated generalization of the family of Riemann surfaces is already known at least for simple situations.

 In Section~\ref{s3.4} we discuss an approach to the problem of constructive description of
$S$-compacta based on the embedding of the set of Robin measures associated with compacta $F_f,$ \ $f \in \mathcal A_e$ into a larger space of probability measures in plane, which we call {\it $e$-critical measures}. Those measures constitute a connected finite dimensional variety and its structure may help to better understand the structure of the discrete set of Robin measures for $F_f,$ \ $f \in \mathcal A_e$ (see \cite{MaRa11}).

Later, in Section~\ref{s4} we also use critical measures to study the problem of stability of convergence in Stahl's theorem under variations of the function $f$ preserving the set of branch points.

\subsubsection{Family of Riemann surfaces \ $\mathcal R_f$,\  $f \in \mathcal A_e$ }\label{3.3.3}
The $S$-property \eqref{S} is essentially equivalent to the fact that the real Green function $g(z)$ of the domain $\Omega = \overline\CC\setminus F_f$ has a harmonic extension to a hyperelliptic Riemann surface ${\mathcal R} = \mathcal R_f$, which may be defined as the Riemann surface of the function $\sqrt{V/A}$ with $V = V_f.$ We interpret $\mathcal R_f$ in a standard way as a two sheeted branched covering over $\overline\CC$. Formula \eqref{Gre} provides a constructive form of this extension.

Recall that on any hyperelliptic Riemann surface $\mathcal R$ there exist a unique function $g = g_{\mathcal R}: \mathcal R \to \mbb R$, which is uniquely defined as harmonic function on the finite part of $\mathcal R$ with asymptotics $g(z) = \log|z| + o(1)$ as $z \to \infty^{(1)}$ and $g(z) = - \log|z| + o(1)$ as $z \to \infty^{(2)}$ and with normalization $g(z^{(1)}) + g(z^{(2)}) \equiv 0$.
We call this function $g = g_{\mathcal R}$ the $g$-function for the Riemann surface $\mathcal R.$ Continuation of the Green function $g(z)$ from the domain $\overline\CC\setminus F_f$ with $f \in \mathcal A_e$ to the Riemann surface $\mathcal R_f$ is exactly the $g$-function for $\mathcal R_f$.

Consequently, the complex Green function $G$ in \eqref{Gre} has
(multivalued) analytic continuation to ${\mathcal R}$, which is a standard
third kind Abelian integral on ${\mathcal R}$ with (logarithmic) poles
at $\infty^{(1)}$ and $\infty^{(2)}$ and divisor $1, -1$ (we call it $G$-function for  ${\mathcal R}$).

Representation $G'(z) = \sqrt{V(z)/A(z)}$, where $V = V_f$, asserted in \eqref{Gre} for $z\in \mbb C \setminus F_f$ is valid for $z \in \mathcal R.$ The extremal compact set $F_f$ is the projection of the zero level $\{\zeta: g(\zeta) = 0 \} \subset \mathcal R$ of $g$-function onto the (extended) plane $\overline \CC$.

Now, together with the collection of $S$-compacta and associated family of polynomials $\widetilde {V}(e)$ we also have the family of Riemann surfaces
$\widetilde{\mathcal R}(e) = \left\{\mathcal R _f : \,\, f \in \mathcal A_e \right\}$.
Next, we will obtain a representation of $\widetilde {V} (e)$  in terms of this family.

Recall that we begin our constructions with a fixed polynomial $A(z) = z^p +\dots$ having simple roots. Next, consider a variable polynomial $V(z) = z^{p-2} +\dots$ (at the moment we do not have relate $V$ to the constructions above). However, assume for now, that the zeros of $V$ are simple and $A$ and $V$ do not have common zeros. Then, Riemann surface of the function $\sqrt{V/A}$ is a generic hyperelliptic Riemann surface of genus $p-1$ having $2p-2$ quadratic branch points at zeros of $AV$. It is well known that the $G$-function for such surface may be written in the form
\begin{equation} \label{Gre-1}
G(z) = G(z; V) = \int_a^z \frac{W(t) dt}{\sqrt{A(t)V(t)}}, \quad \text{where}\quad W(z) =  z^{p-2} +\dots
\end{equation}
($a$ is a root of $A$).  The polynomial $W$ is uniquely determined by the polynomial $AV.$ Since $A$ is fixed, this defines a mapping $\Phi: V \to W$. It is yet defined under the assumption that zeros of $V$ are simple and different from roots of $A$, but the mapping has continuous extension to the whole space $ \mbb P_{p-2}^{(1)}$ of monic polynomials of degree $p-2$ (actually, we need only the restriction of the mapping $\Phi: \mbb P_{p-2}^{(1)} \to  \mbb P_{p-2}^{(1)}$ to the space of polynomials with zeros in the convex hull of roots of $A$).

It follows from \eqref{Gre-1} that $G' = W/\sqrt{AV}$. On the other hand (we return to original settings), if  $V \in  \widetilde {V} (e)$ then we have $G'(z) = \sqrt{V(z)/A(z)}$ according to  \eqref{Gre}. Combining the two representation we obtain $W = V$. In other words, $V\in\widetilde {V} (e)$
implies that the polynomial $V$ is a fixed point of the mapping
$\Phi$. Reciprocally, any fixed point of $\Phi$ is in $ \widetilde {V} (e)$
and, therefore, $ \widetilde {V} (e)$ is equivalently defined as set of fixed
points of the mapping $\Phi$.

It is generally possible that polynomial $V$ has common zeros with $A$.
Then those common zeros are canceled in the ratio $V(z)/A(z)$ and the
problem is reduced to a similar one with a smaller set $e$  (of roots
of $A$). Such reduction would not be a significant event. For instance,
let $p=3$ and roots of $A$ be collinear. Then the root of $V$ will
cancel the middle root of $A$ and the problem reduces to the one with
$p=2$. Cancellations of the other two roots of $A$ are banned by
assumption that all roots of $A$ are branch points of the functions $f
\in \mathcal A_e.$

Reduction of the genus of the surface $\mathcal R _f$ may also be the result of the presence of multiple zeros of $V$ and this is a common occasion, which has an important meaning. All, except for maybe one, polynomials $V\in \widetilde {V} (e)$ have multiple roots. Loosely speaking, this fact is a reflection of a possible variety of branch types of the functions $f \in \mathcal A_e.$ Anyway, combinatorics of the set $\widetilde {V} (e)$ is in part determined by multiple roots of $V$.

Suppose that $f$ has a generic branch type, that is, continuation along any nontrivial loop leads to a different branch. Then, associated $S$-compact $F_f$ is a continuum; it is, therefore, the Chebotarev continuum for $e$.  In a situation of a ``common position'' for configuration of points in set $e$, polynomial $V^0$  associated with the Chebotarev continuum will have simple zeros (it may be viewed as a definition of a ``common position''). In such a situation $V^0$ is the only fixed point of the mapping $\Phi$ with simple zeros. All other polynomials $V \in  \widetilde {V} (e)$ will necessarily have multiple roots and, so, reduced genus of associated Riemann surfaces. We do not go into further details. Discussion of the structure of set $\widetilde {V} (e)$ is continued in the next section from a different point of view.

In conclusion of this section we make the following remark. In case of a finite set $e$ the Stahl's theorem may be equivalently formulated in terms of convergence of the remainder $R_n$ in \eqref{eq2.2} on a Riemann surface (in particular, this gives an alternative approach to the way of introduction of the Riemann surface $\mathcal R = \mathcal R_f$).

The theorem may be stated as follows. For a given $f\in \mathcal A_e$ there exist a hyperelliptic Riemann surface  $\mathcal R$ such that (with proper normalization) sequence $\frac 1n \log |R_n|$ converges in capacity on $\mathcal R$ to the $g$-function of $\mathcal R$. In the equivalent form: the sequence of normalized logarithmic derivatives $\frac {R'_n }{n R_n}$ converges to $G'(z)$ in the plane measure on $\mathcal R$. Then, the surface $\mathcal R$ is  uniquely defined by the additional conditions that the projection $F$ of zero level of $g$ onto the plane makes $f$ single-valued and also a jump of $f$ across any arc from $F$ is not identical zero.

It is possible that in such form the Stahl's theorem may be directly generalized for the first kind Hermite--Pad\'e approximants for systems of functions with branch points.

\subsection{Critical measures \texorpdfstring{$\mc M_e$}{Me}}\label{s3.4}
For a finite set $e=\left\{a_1,\dotsc,a_p\right\}$ we define
$e$-critical measures as critical points of the energy functional
$\mc{E}(\mu) = - \int \log|x-y|d\mu(x) d\mu(y)$ with respect to local
variations with fixed set $e$. More exactly, for a smooth complex
function $h(z)$ in a neighborhood of $\supp\mu$ we define point
variations $z\to z^t=z+th(z)$, where $|t|\in(0, \epsilon)$, and then
variations of the measures $\mu \to \mu^t$ by $d\mu(z) = d\mu^t(z^t)$.

An associated variation of energy (derivative in the direction $h$) is defined by 
\begin{equation}
\label{Var}
D_h\mc{E}(\mu)=\lim_{t\to 0+}\frac{1}{t}\(\mc{E}\(\mu^t\)-\mc{E}(\mu)\).
\end{equation}
Finally, we say that $\mu$ is $e$-critical, if for any function $ h $ satisfying condition $h(a) = 0$ for any $ a\in e$,
we have $ D_h\mc{E}(\mu)=0$.
The set of all such measures is denoted by $\mc M_e$.

Critical (stationary) measures were first introduced in~\cite{GoRa87}
and then used in~\cite{PeRa94}. A systematic study of critical measures
(with rational external fields) was presented in~\cite{MaRa11} in
connection with zero distribution of Heine--Stieltjes polynomials; see
also review~\cite{Rak12}.
Here we use critical measures as an approach to describe the set of
Robin measures of $S$-compact sets $F_f$ associated with a fixed set $e$. Later in
Section~\ref{s4} they are also used to study stability of convergence in Stahl's theorem.

It is important to observe, first, that Robin measures of all $S$-compacta $F_f$ are $e$-critical measures; second, basic properties of Robin measures of $S$-compacta are preserved for critical measures.
In particular, the potential of any $e$-critical measure $\mu$ has the $S$-property
presented by \eqref{S} with $U^\mu$ in place of $g$.
Next, for any critical measure $\mu$ there exist a polynomial
$\ V(z) = \prod_{j=1}^{p-2} (z-v_j)$ such that with $\ A(z) = \prod_{k=1}^{p} (z-a_k)$ we have
\begin{equation}
\label{eq4.6}
U^\mu(z)=\Re\int_{a_1}^{z}
\sqrt{{V(t)}/{A(t)}}\,dt,\quad
d\mu(z)=\frac{1}{\pi}\left|\sqrt{{V}/{A}}\,dz\right|.
\end{equation}
Moreover, $\supp\mu$ is a union of critical trajectories of
$-\({V(z)}/{A(z)}\)(dz)^2$ and this differential has closed
trajectories just as for the Robin measures of $S$-compact sets $F_f$.
Finally, both sets of measures may be characterized in terms of the associated polynomials $V$.

Using the zeros $v_j$ of $V$ as parameters we represent the set $\mc M_e$ of critical measures 
as a subset in the space of vectors $\{ v = \(v_1,\dotsc,v_{p-2}\)\}$
from $\mbb C^{p-2}$. In these coordinates $\mc M_e$ is represented as a union
of $3^{p-2}$ bounded bordered domains, which we call cells. Each cell is
a bounded bordered manifold of real dimension $p-2.$  Interior points of each cell
correspond to measures $\mu$, whose support $\Gamma=\supp\mu$ consists of
exactly $p-2$ simple disjoint analytic arcs $\Gamma_j$ with endpoints
from the set $\{a_k,\ v_j\}$. Finally, $v$ coordinates
of measures $\mu \in\mc M_e$ are defined by systems of equations
\begin{equation}
\label{Equ}
\Re\int_{\Gamma_j}\sqrt{{V(t)}/{A(t)}}\,dt=0,\  j=1,\dotsc,p-2; \quad V(t) = (t-v_1) \cdots (t-v_{p-2}).
\end{equation}
Particular cell is identified by homotopic type of arcs $\Gamma_j$.

Robin measures of $S$-compact sets are among $e$-critical measures and
their representations in terms of $v$-coordinates are located on boundaries
of cells.

The space $\mc M_e$ is connected and each critical measure may be, in a
standard way, connected with the Chebotarev continuum associated with
$e$, which may be defined as the only continuum (closed connected set)
in the set of $S$-compacta for functions $f\in \mc A_e$. The roots of the polynomial
$V_0 \in\widetilde{V}(e)$, associated with the Chebotarev continuum, effectively
play the role of origin in the $v$-coordinate system and
the corresponding ``deformation theory'' is in part described in~\cite{MaRa11}.

Further, equilibrium measures of $S$-compact sets satisfy \eqref{Equ}
and also $p-2$ additional equations, which distinguish them among all
critical measures. Potential of any interior critical measure $\mu$,
supported on arcs $\Gamma_j$, keeps constant value on those arcs, that
is, we have
$$
U^\mu(z)=C_j,\quad z\in\Gamma_j, \quad \supp\mu=\bigcup\limits_{j=1}^{p-1}\Gamma_j.
$$
The collection of constants $C = \(C_1,\dotsc,C_{p-2}\)$ may be used to
parameterize points in a particular cell in $\mc M_e$ (note that only
$p-2$ of the constants are independent). The differences of those
constants correspond to the parametrization of a cell by ``height of
cylinders'' in terms of the general moduli problem, which may be
associated with critical measures. There is a dual parametrization by
``lengths of circles'' which corresponds to masses $\mu(\Gamma_j)$ (see
~\cite{Str84}).

Now, additional equations, which determine equilibrium
measures of $S$-compacta in terms of $C$-coordinates, are
$$
C_1=\dotsb=C_{p-1}.
$$
We have a total of $2p-2$ real equations for the same number of real
parameters in $V$. For some further details
see~\cite{Str84},~\cite{MaRa11},~\cite{Rak12}.

\section{Some generalizations and conjectures}\label{s4}

The method outlined above may be developed in several directions. In this
section we make a few remarks related to possible generalizations.

\subsection{Dependence of Pad\'e denominators from the function}\label{s3.1}
Let $f\in \mc A_e$ be a function with a finite set $e$ of branch points
defined by an element \eqref{f} at infinity. Let $\pi_n(z)=(P_n/Q_n)(z)$
be the associated Pad\'e approximants at infinity.

Suppose that we make a small variation of the function $f$ in class $\mc A_e$.
In other words, consider a  function $\widetilde f \in \mc A_e$,
which is close to $f$ in some sense, the location of its branch points is the same,
but their character may change.

We want to figure out how much the denominator $Q_n = Q_n(\widetilde f)$ will change,
say, for a fixed large enough $n$. To be more precise, here we have in mind a significant change,
and as a first step toward investigation of the problem, we will discuss a possible change in the
limit zero distribution.  Since the rate of convergence in Stahl's theorem is determined by
the limit zero distribution, the problem is essentially about the (rough) stability of the convergence
in this theorem.

It turns out that the answer depends on what exactly was the $S$-compact set $F_f$
for the function $f$. As usual, assume for simplicity that all points $a_k\in
e$, $k =1,\dots, p$ are actual branch points for $f$.

If $F_f$ was the Chebotarev continuum $F_e$ for $e$ (in other words,
function $f$ has a ``generic branch type'') then small variations of $f$
will not produce a dramatic effect. It is not difficult to see that any such small
enough variation of $f$ will remain to be of a ``generic branch type'' and,
therefore, will have the same $S$-compact and the same limit distribution
(the Robin measure of this compact). Thus, the dependence of $Q_n$ on $f$ is
essentially continuous (and asymptotically continuous).

If $F_f$ was any other $S$-compact set, then this dependence
is not continuous, since dependence $F_f$ from $f$ is not continuous if $F_f$ was not the Chebotarev compact for $e$.
The branch types, which may be obtained by small variations of $f$, depend on
$F_f$. Anyway, it is clear that the generic branch type may be obtained form any other one using arbitrary small variations and it is enough to prove discontinuity of $F_f$ as function of $f$ at any $f$ whose $S$ compact is not the Chebotarev's one. All facts above are still corollaries of Stahl's theorem.

The situation changes if we consider a sequence of variations depending on $n$, which converges to zero
as $n\to\infty$. What exactly will happen with zero distribution depends on characteristics of the function, characteristics of the variation and the relation between $n$ and the magnitude of variation (it is possible to consider also variations of the locations of branch points, but the effect will be similar).

Formally, let   $f,\ f_n\in \mc A_e$ for $n\in \mbb N$ and the sequence $f_n$ converges to $f$ as $n\to\infty$.
Let $Q_n = Q_n(\widetilde f_n)$. From what was said above follows that, if $F_f$ is the Chebotarev continuum for $e$, then the sequence $\frac1{n}\mc{X}\(Q_n\)$ converges weakly to the Robin measure of this continuum.

If $F_f$ is  different from the Chebotarev continuum for $e$ the sequence $ \ \frac1{n}\mc{X}\(Q_n\)$ is not generally weakly convergent. We may claim  that only the weak-$*$ limit of any convergent subsequence belongs to the set $\mc M_e$ of critical measures for $e$. The measures, which are included in the limit set for a given $f$, depend on $ F_f$ and the character of convergence.
However, using different functions $f\in \mc A_e$ we may
obtain any $\mu \in \mathcal M_e$  as a limit along the whole sequence.
In other words, any $\mu \in \mathcal M_e$ is a weak limit of the whole sequence
$ \ \frac1{n}\mc{X}\(Q_n\)$ for some selection of functions $f,\ f_n\in \mc A_e$.

To state a theorem formally presenting assertions above we have to define the convergence $f_n\to f.$
 We give a simple example of such a theorem with particularly simple kind of convergence.
Consider the following model class of functions
\begin{equation}
\label{Cls}
\mc L_e = \left \{f:f(z) = \prod_{k=1}^p (z - a_k)^{\alpha_k} \right\}, \qquad \sum_{k=1}^p \alpha_k = 0.
\end{equation}
We assume that $e = \{a_k\}$ is fixed and $\alpha_k$ are parameters; as
usual we assume that each $a_k$ is an actual branch point of the
function ($\alpha_k$ is not an integer). We have $\mc L_e \subset \mc
A_e$ and the class $\mc L_e$ is representative enough in the sense that all
possible branch types are presented by functions from $\mc L_e$.
The convergence  $f_n \to f$ for functions from $\mc L_e$ is understood as the convergence
$\alpha_{k,n} \to \alpha_k$ of the $\alpha$-parameters of $f_n$ to those of $f$.

Now we can state a version of the theorem related to the class $\mc
L_e$.

\begin{thm}\label{Stab}
Let the sequence of functions $f_n \in \mc L_e$ converge to $f\in \mc L_e$.

If the extremal compact set $F = F_f$ for $f$ is the Chebotarev
continuum $F_e$ for $e$, then $\frac1{n}\mc{X}\(Q_n(f_n)\)\ \
{\overset{*}{\to}}\ \ \lambda$, where $\lambda$ is the Robin measure
for~$F$.

For any $\mu\in \mc M_e$ there exist a convergent sequence of functions $f_n \in \mc L_e \to f \in \mc L_e$ such that
$\frac1{n}\mc{X}\(Q_n(f_n)\)\ \ {\overset{*}{\to}}\ \ \mu$ as $n\to\infty$.
\end{thm}

Theorem~\ref{Stab} generalizes the Stahl theorem in
the same way that Theorem~1 in~\cite{GoRa87} generalizes the
$\rho^2$-theorem. Proofs of both theorems may be based on the description
of the set of critical measures $\mc M_e$ outlined above and also on
Theorem~1 from~\cite{GoRa87}. In the next sections we briefly discuss
this theorem and some of its possible generalizations. At the same time Theorem~\ref{Stab} may be proved in a very simple way using Laguerre type differential equation for Pad\'e denominators of functions
from $\mc L_e.$

\subsection{Conjectures on zero distribution of complex orthogonal polynomials}\label{s4.2}
Here we present some conjectures connected to Hermite--Pad\'e polynomials (their circle of applications may be larger). Thus, we touch a general problem of generalizations of the theory outlined above in this paper for the case of Hermite--Pad\'e polynomials. This is one of the central problems in the theory and at the moment the problem is essentially open.

As a starting point we need a version of a general  theorem from~\cite{GoRa87} (GRS theorem),
which is for the moment, probably, the most advanced known theorem related to zero distribution of complex orthogonal polynomials.

\subsubsection{GRS theorem}\label{s4.2.1}
To state the theorem we need the following definition.

We say that  a compact $F\subset \mbb C$ has $S$-property in an external field $\varphi$ harmonic
in a neighborhood of $F$, if equality in \eqref{S} holds for $g = U^\lambda + \varphi$ -- total
potential of equilibrium measure $\lambda=\lambda_{\varphi,F}$ for $F$ in
the external field $\varphi$. The $S$-property implies that $F$ is at most
countable union of disjoint open analytic arcs $F^0$ and a set of capacity
zero (here we assume from the beginning that an $S$-compact associated with the problem exists).

Now, we state assumptions of the theorem.

We assume that we are given a domain $\Omega$ in $\mbb C$, a compact set $F$ in $\Omega$
and a sequence of functions $\Phi_n(z)\in H(\Omega)$, which
converge $\Phi_n(z)\to\Phi(z)$ uniformly on compact subsets of $\Omega$ as
$n\to\infty$.

Assume that $F$ has $S$-property in the external field $\varphi=\Re\Phi(z).$

Further, let $f\in H(\Omega\sim F)$ be a function whose jump across any
arc from $F^0$ is not identical zero and polynomials $Q_n(z) \in \PP_n$ be defined by orthogonality relations with weights $f_n=f e^{-2n\Phi_n}$
\begin{equation}
\label{Ort-G}
\oint_{F} Q_n(z)\,P(z)\,f_n(z)\,dz=0,\qquad
\text{for any polynomial}\quad P \in \PP_{n-1}.
\end{equation}
Integration in \eqref{Ort-G} goes along a contour(s) in $\Omega\setminus F$, homotopic to the boundary of $\ol{\mbb C}\setminus F$.

Finally, assume that the complement to the support of the equilibrium measure $\lambda=\lambda_{\varphi,F}$ for $F$ in the external field $\varphi$ is connected.

The following is Theorem~1 from~\cite{GoRa87}.

\begin{thm}\label{thm3}
Under the assumptions above we have \ $\frac1n\,\mc{X}\(Q_n\)\overset{*}{\ \to\ }\lambda$.
\end{thm}

Orthogonality conditions in \eqref{Ort-G} above are rather general, but in a number of situations the theorem may not be directly applied.  It happens often in the study of zero distribution of the Hermite--Pad\'e polynomials.  These polynomials are defined by systems of orthogonality relations and reduction of such systems to orthogonality with respect to single weight (if possible) lead to more general forms of orthogonality. Next we give two comparatively simple examples of different nature.

\subsubsection{A conjecture related to Hermite Pad\'e polynomials for a Nikishin system}\label{s4.2.2}
In many cases study of Nikishin systems may be reduced to a problem of asymptotics for orthogonal polynomials $Q_n$, defined by relations similar to \eqref{Ort-G} in theorem \eqref{thm3}, but with the weight functions $f_n$ in \eqref{Ort-G} depending not only on $n$ but also on
 the polynomial $P$. In other words, polynomials $Q_n$ are orthogonal to
 some collection of functions, which are not pure polynomials, but
polynomials $P$ with multiplier $f_n$, depending not only on $n$, but also on $P$.

We formulate a conjecture for the case when only $\Phi_n$ depend on $P$.

Let polynomials $Q_n$ satisfy orthogonality conditions \eqref{Ort-G}
 with $\Phi_n(z) = \Phi_n(z; P)$.
All the assumptions of Theorem~\ref{thm3} above are preserved.
In addition, we assume that for any
sequence of polynomials $P \in \PP_{n-1}$ such that $\frac1{n}\
\mc{X}\(P_n\)\ \ {\overset{*}{\to}}\ \ \lambda$ we have $\Phi_n(z; P_n))
\to \Phi(z)$.

\begin{conj}\label{conj1}
Under the above assumptions we have
$\ \frac1n\mc{X}(Q_n)\overset{*}\to\lambda$.
\end{conj}

Conjecture~\ref{conj1} is a part of joint work with S.~Suetin (in
progress). It is partially suggested by the results of the
paper~\cite{RaSu13}, where Hermite--Pad\'e polynomials of the first
kind were considered for a Nikishin system of two Markov-type functions
$f_1, f_2$ on the union $E$ of a finite number of disjoint real closed
intervals $E_j.$ We outline settings of the paper without going into
all the details related to the situation.

In the paper~\cite{RaSu13} it was assumed that the ratio of two jumps
$f(x):=\Delta f_2(x)/\Delta f_1(x)$, $x\in{E}$, is an analytic complex-valued
function on $E$ and $f$ has an analytic continuation from each $E_j$ along any
path in $\myo\CC$ avoiding the finite set $e_f$ of the branch points of $f$.
It was also assumed, that the set $e_f$ is symmetric with respect to
real axis.

First it was proven that (under some additional technical assumptions)
there exists a unique compact set $F$, such that
$f\in H(\myo\CC\setminus{F})$ and $F$ has $S$-property with respect
to some related equilibrium problem for a mixed Green-logarithmic potential.

Let $Q_{n,0},Q_{n,1},Q_{n,2}\in\PP_n$, $Q_{n,2}\not\equiv0$, be the
Hermite--Pad\'e polynomials of the first kind for the system $[1,f_1,f_2]$, that is,
the following relation holds
\begin{equation}
\label{HP-I}
(Q_{n,0}\cdot 1+Q_{n,1}f_1+Q_{n,2}f_2)(z)=O\(\frac1{z^{2n+2}}\),\quad z\to\infty.
\end{equation}
The following orthogonality relation of type~\eqref{Ort-G} for the polynomial $Q_{n,2}$ was obtained
in~\cite[formula~(119)]{RaSu13},
\begin{equation}
\oint_{F}Q_{n,2}(z)P_n(z)
\biggl\{h_{n+m}(z)\frac{\tau_n^2(z)}{q_{n}(z)}
\int_E\frac{q_{n}^2(\zeta)\tau_n^2(\zeta)}{z-\zeta}\frac{d{\mathfrak m}_n(\zeta)}{P_n(\zeta)} f(z)\,
\biggr\}\,dz=0
\label{RaSu_Ort-G}
\end{equation}
($P_n\in\PP_{n-1}$ is an arbitrary polynomial). Finally, these orthogonality relations were used to prove that the sequence $\frac1n\mc{X}(Q_{n,2})$ weakly converges to the equilibrium measure for the problem mentioned above; for more details see~\cite{RaSu13}.

Connection of this result with the Conjecture \ref{conj1} is established by the following fact. The function in curly brackets in \eqref{RaSu_Ort-G}, which plays role of multiplier for $P_n$, satisfies conditions in Conjecture \ref{conj1}. Thus, the theorem above supports the conjecture.

\subsubsection{A conjecture on incomplete complex orthogonal polynomials}\label{s4.2.3}
Hermite--Pad\'e polynomials also lead to another type of asymptotics problems for orthogonal polynomials. Before discussing this problem (in the last section below), we introduce an auxiliary problem for complex orthogonal polynomials. The problem may have, however, an independent value. We restrict ourselves with the simplest possible version of the problem.

Let $ f(z) \in \mc A_e$, where $e=\{a,b\}$. That is, function $f$
(defined by an element at infinity) has two branch points at $a$ and $b
\ne a$. Let $N$ and $n \leq N$ be two natural numbers and the two polynomials
$Q_N \in \mbb P_N$ satisfy relations
\begin{equation}
\label{Ort-L}
\oint_F Q_N(z)\,P(z) \,f(z) \,dz=0,\qquad \text{for any polynomial}\quad P \in \mbb P_{n-1},
\end{equation}
where $F$ is a curve connecting $a$ and $b$. Note that here we do not
assume that a special curve is given.  Any curve  $F$ from the class $\mc F$ of curves connecting   $a$ and $b$ may be used in \eqref{Ort-L} by the Cauchy integral theorem. Finding a special curve  $\Gamma \in \mc F$ will be a part of the problem.

Suppose that $\, n, N \to \infty$ in such a way that $N/n \to k
>1$. What can be said about the zero distribution of $Q_N$?

Clearly, under these assumptions the polynomial $Q_N$ is not uniquely defined
and we cannot expect that the sequence of counting measures
$\frac1n\,\mc{X}\(Q_N\)$ is convergent. Instead, we suggest that any limit
point of this sequence satisfies certain inequality. To state formally this
inequality, we first need to select a convergent subsequence
\begin{equation}
\label{Conv}
\frac1n\,\mc{X}\(Q_N\)\overset{*}{\ \to\ }\mu \qquad \text{as}\quad n\to \infty, \quad n\in \Lambda
\end{equation}
($\Lambda$ is a sequence of natural numbers). Since $N/n\to k=1$, we have $\mu
(\mbb C) = k >1.$

The potential $\varphi = U^\mu$ of $\mu$
will play the role of external field in the problem we are going to
consider. We denote by
$$
\mc E_\mu(\nu)  = \mc E(\nu) + 2 \int U^\mu\, d \nu
$$
the weighted energy of a measure $\nu$ in the external field $\varphi.$
Note that here and in the sequel we use the abbreviated notation: $\mc
E_\mu(\nu)$ stands for $\mc E_\varphi (\nu)$ with $\varphi = U^\mu$ (compare to \eqref{EqM}).

For a fixed $F \in \mc F$ we denote by  $\lambda_{F,\, \mu} \in \mc {M}(F)$
the minimizing (equilibrium) measure on $F$ in the external field $\varphi
= U^\mu$ and the equilibrium energy by
$$
  \mc E_ \mu (\lambda_{F, \,\mu})
= \min_{\nu\in\mc {M}(F)}\ \mc E_\mu(\nu) ,
$$
where $\mc {M}(F)$ is the set of probability measures on $F$.

Next, we introduce the functional of equilibrium energy $\mc E_\mu [F]$ and assert existence
of a compact set $\Gamma = \Gamma_\mu  \in \mc F$ maximizing this functional (see \cite{Rak12})
\begin{equation}
\label{Ext}
\mc E_\mu \, [\Gamma] \,  = \, \max_{F\in \mc F}\, \mc E_\mu [F] \qquad \text{where} \qquad   \mc E_\mu [F ] =  \mc E_ \mu (\lambda_{F, \,\mu}).
\end{equation}
Finally, we define a mapping $\mu \to \lambda$ in the space of probability measures in the plane by
\begin{equation}
\label{Func}
\lambda(\mu) = \lambda(\mu, \mc F) = \lambda_{\Gamma, \,\mu},
\end{equation}
where $\Gamma = \Gamma_\mu$ is the extremal compact in \eqref{Ext}. The conjecture is formulated in terms of this function.

\begin{conj}\label{conj2}
For any subsequential limit $\mu$ of the sequence $\ \frac1n\,\mc{X}\(Q_N\)$
we have $ \mu \geq \lambda(\mu)$.
\end{conj}

In a number of situations, Conjecture \ref{conj2} can be proven. We mention one
such situation where the proof may be obtained using the GRS-method. Suppose
that the limit distribution is known for a part of zeros containing
$N - n$ zeros. Let this limit distribution is represented by a known measure $\sigma$.
In other words, we assume that $N - n$ have factorization $Q_N = q_n g_n$, where
the sequence of polynomial $g_n \in\PP_{N-n}$ has a limit distribution $\sigma$;
formally $\frac1n\,\mc{X}\(g_n\)\overset{*}{\ \to\ }\sigma. $

Suppose also, that the class $\mc F$ of continua $F$ connecting $a$ and
$b$ contains a continuum $\Gamma$ with $S$-property in the external field
$\varphi(z) = \frac 12 U^\sigma(z)$.  Then the sequence
$\frac1n\,\mc{X}\(q_n\)$ is weakly convergent to $\lambda =
\lambda_{\varphi,\,\Gamma}$ according to Theorem~\ref{thm3}. It follows
that the sequence $\frac1n\,\mc{X}\(Q_N\)$ converges to $\mu = \lambda +
\sigma$, and finally the sequence $\frac1n\,\mc{X}\(Q_N\)$  converges to $\mu
\geq \lambda$. Since the equilibrium measure of an $S$-compact has the
$\max$--$\min$-property, we have $\lambda = \lambda(\mu)$ and the assertion of
Conjecture~\ref{conj2} follows.

Thus, in this situations the $\max$--$\min$ definition of $\lambda(\mu)$ can
be equivalently formulated in terms of the $S$-property. In general, we have to
define $\lambda(\mu)$ in terms of ``$\max$--$\min$'', since the external fields
associated with the problem may not be harmonic (even smooth) around the
extremal compact.

\subsubsection{A conjecture related to Hermite--Pad\'e polynomials for an Angelesco system}\label{s4.2.4}
As an example of possible application of Conjecture \ref{conj2} we mention the
problem of zero distribution for denominators of the second kind
Hermite--Pad\'e approximants in Angelesco case.

The Simplest settings are as follows. Let $e_1 =\{a_1, b_1\}$
and $e_2 =\{a_2, b_2\}$ be two sets, where $a_i \ne b_i$ for $i = 1,2$ are given. Then, two
functions $f_1 \in \mathcal A_{e_1}$ and $f_2 \in \mathcal A_{e_2}$ are
defined by their Laurent series at infinity. Assume that $\{a_i, b_i\}$ are
actual branch points of $f_i$.  Finally, a nontrivial polynomial $Q = Q_{2n}
\in\PP_{2n}$ is defined by the pair of conditions
$$
(Q f_1 - P_1)(z)  = O\left(z^{n+1}\right), \qquad (Q f_2 - P_2)(z) = O\left(z^{n+1}\right)
$$
as $z \to \infty$, where $P_i$ is the polynomial part of $Qf_i$ at infinity
($i = 1,2$).

We assume that the couple of functions $f_1, f_2$ (or, rather, couple of
sets $e_1, e_2$) present ``Angelesco case'', which  informally speaking means
that $e_1$ and $e_2$ are ``well separated'' (far enough from each other).
The formal definition is presented below after related definitions are
introduced. As an example, we note that if all branch points are real, then we
define Angelesco case by the condition that the intervals $(a_1, b_1)$  and
$(a_2, b_2)$ are disjoint. It is known that in such case the limit zero
distribution of the sequence $Q_{2n}$ is defined by a matrix equilibrium
problem on a pair of conductors $F_1 = [a_1, b_1]$  and $F_2 = [a_2, b_2]$
(see \cite{GoRa81} and \cite{GoRaSo97} for Markov case). For the complex case we
have to use a matrix $S$-equilibrium problem, which is defined below.

For $i = 1,2$ denote by $\mc{F}_i$ the class of continua in the plane connecting points  $a_i$ and $b_i$.
We consider the class of vector compacts $\vec{\mc F}=\(F_1, F_2\)$, where $F_i \in \mc{F}_i.$
For a fixed vector-compact set $\vec{F}=\(F_1, F_2\)\in\mc{F}$ we define the class of vector-measures
$$
\vec{\mc{M}}
=\left\{ \(\mu_1, \mu_2\): \, \mu_j \in\mc{M}\(F_j\)\right\},
$$
where $\mc{M}\(F_i\)$ is the set of probability measures on $F_i$.
The energy of the vector measure $\vec{\mu} = \(\mu_1, \mu_2\)$ is defined by
$$
\mc{E}\(\vec{\mu}\)
=[\mu_1,\mu_1] + [\mu_1,\mu_2] + [\mu_2,\mu_2] ,\qquad \vec{\mu} = \(\mu_1, \mu_2\),
$$
where $[\mu,\nu]=\int V^\nu d\mu$ is the mutual energy of $\mu$ and $\nu$.
In a more general situation, the energy of a vector measure is defined by a
matrix $A$ with constant elements $a_{ij}, \, i,j = 1,2$, so that the
matrix-energy is  $\mc{E}\(\vec{\mu}\) = \sum a_{ij} [\mu_i, \mu_j]$. In
our case the elements of matrix $A$ are $a_{11}= a_{22} = 1$  and $a_{12}=
a_{21} = 1/2$. This is the positive definite matrix and, moreover,  $a_{ij}\ge
0$. It follows that for any $\vec{F}=\(F_1, F_2\)\in\mc{F}$ there exists a
unique $\vec{\lambda}\in\vec{\mc{M}}$, such that
$$
\mc{E}[\vec{F}] = \mc{E}(\vec{\lambda})=\min_{\vec{\mu}\in\vec{\mc{M}}(\vec F)}\,\,\mc{E}(\vec{\mu}),
\qquad \vec{\lambda} = \(\lambda_1, \lambda_2\).
$$
The vector-measure $\vec{\lambda}$ is the equilibrium measure for $\vec{F}$
associated with matrix $A$; \ $\mc{E}[\vec{F}] $ is the equilibrium energy
of $\vec F$ (see original papers \cite{GoRa81},~\cite{GoRa85},~\cite{GoRaSo97})
and recent developments in  \cite{Bec13},~\cite{HaKu13}).

Further, there exists a vector-compact set $\vec{\Gamma}=\(\Gamma_1,
\Gamma_2\)\in\mc{F}$ maximizing the equilibrium energy
$$
 \mc{E}[\vec{\Gamma}] = \max_{\vec{\Gamma} \in \vec{F}} \,\, \mc{E}[\vec{F}].
$$
The existence of maximizing vector-compact sets $\vec{\Gamma}$ may be proved by
the methods presented in \cite{Rak12}. In general, it is not unique, but the
associated equilibrium measure $\vec{\lambda} = (\lambda_1, \lambda_2) $ is
unique.

What can be asserted about the limit zero distribution of Hermite--Pad\'e
denominators $Q_{2n}$ essentially depends on the structure of
$\vec{\Gamma}$ or, better to say, the structure of $\vec{\lambda}$. If the supports
of  $\lambda_1$ and $\lambda_2$ are essentially overlapping, then the vector
measure $\vec{\lambda}$ does not describe the zero distribution of the polynomials
$Q_{2n}$ and the case under consideration is not an Angelesco case. In such
situation the equilibrium problem has to be modified; we refer to papers
\cite{Apt08} and \cite{ApKuAs08} for further details.

If the supports of $\lambda_1$ and $\lambda_2$ are not intersecting, then we
have Angelesco case and we assume this condition in what follows (the case when
there is a small -- say, finite -- intersection may be included, but we
restrict our considerations to the disjoint situation). Now, the main
hypothesis on the zero distribution of Angelesco Hermite--Pad\'e polynomials
is stated as follows.

\begin{conj}\label{conj3}
We have \ $\frac 1n \mc{X}(Q_{2n}) \ \ {\overset{*}{\to}}\ \  \lambda_1
+ \lambda_2$ \ where  \ $\vec{\lambda} = (\lambda_1,  \lambda_2)$
\ is the equilibrium measure of the extremal compact
$\vec{\Gamma}=\(\Gamma_1, \Gamma_2\)$.
\end{conj}

In a number of cases the conjecture is proved under some additional restrictions. First, if the sets $e_1$ and $e_2$ are far enough
the proof may be given on the basis of the GRS-method. This approach is simple and general. The number of functions and the number and character of the branch points are actually not important if the sets $e_k$ of branch points are well separated. However, it is difficult
to obtain sharp estimates of the critical distance between sets.

Second, with some additional assumptions on the character of the branch points
the strong asymptotics for $Q_{2n}$ has been derived in \cite{ApKuAs08} for
two functions, each with two algebraic-logarithmic branch points. The proof
uses the steepest descent method for matrix Riemann--Hilbert
representation of $Q_{2n}.$  This method is sensitive to the number of functions and the
numbers of their branch points.  It is not clear if it may be generalized
for arbitrary branch points.

Now, a proof for an arbitrary Angelesco situation may be reduced to
Conjecture~\ref{conj2} (or its proper generalization for more than two functions with any number
of branch points). Such reduction would require some additional
potential-theoretic considerations.

We will mention the shortest way of reduction which is based on one more conjecture.

Let \ $\vec{\lambda} = (\lambda_1,  \lambda_2)$ \ be the equilibrium
measure of the extremal compact set $\vec{\Gamma}=\(\Gamma_1, \Gamma_2\).$ We define
$\mu = \lambda_1+  \lambda_2$ as in Conjecture~\ref{conj3}. Assuming that the Angelesco
case is in effect, the extremal vector-compact set $\vec \Gamma$ has the following
important property: both components $\Gamma_1$ and $\Gamma_2$ have the $S$-property
in the external field $\varphi(z) = U^\mu(z)$. Further, the $S$-property may be
rewritten as the ``energy $\max$--$\min$ property'' and, therefore, the three
measures above satisfy the following relation
\begin{equation}
\label{VecEq}
 \mu \geq \lambda_1 = \lambda(\mu, \mc F_1),\qquad  \mu \geq \lambda_2 = \lambda(\mu, \mc F_2),
\end{equation}
where both measures $ \lambda(\mu, \mc F_1)$ and  $\lambda(\mu, \mc F_2)$
are defined in \eqref{Func} above and $\mc F_i = \mc F(f_i)$ are the classes of
admissible cuts for functions $f_1, f_2$.

Of course, we actually have equality in \eqref{VecEq} but we need
inequalities to make stronger the inverse assertion.

\begin{conj}\label{conj4}
In the Angelesco case, the measure $\mu = \lambda_1 +  \lambda_2$ is the only
positive Borel measure in the plane satisfying \eqref{VecEq} with $\mu(\mbb C)
= 2.$
\end{conj}

If both conjectures~\ref{conj2} and~\ref{conj4} are true, then
Conjecture~\ref{conj3} is also true since it is a direct corollary of the
first two conjectures.

The author thanks the referee for the valuable remarks directed towards improving the style of the paper and also for indicating instances where additional explanations are desirable.


\end{document}